\documentclass[twoside,11pt]{article}
\title{} \author{} \date{}
\usepackage{latexsym,amssymb,times}
\newtheorem{te}{Theorem}[section]

\newtheorem{fac}[te]{Fact}
\newtheorem{lem}[te]{Lemma}
\newtheorem{cla}{Claim}

\def\dok{\noindent{\bf Proof. }}
\def\kdok{\hfill $\Box$ \par \vspace*{2mm} }
\def\a{\alpha}

\def\o{\omega}
\def\e{\eta}
\def\f{\varphi}
\def\k{\kappa}
\def\p{\psi}

\def\t{\tau}
\def\th{\theta}
\def\L{\Lambda}
\def\F{\Phi}
\def\N{\mathbb N}
\def\B{\mathbb B}
\def\P{\mathbb P}
\def\X{\mathbb X}

\def\Q{\mathbb Q}
\def\S{\mathbb S}
\def\A{{\mathcal A}}
\def\D{{\mathcal D}}

\def\CL{{\mathcal L}}
\def\T{{\mathcal T}}
\def\la{\langle}
\def\ra{\rangle}
\def\dom{\mathop{\rm dom}\nolimits}

\def\Bt{\mathop{\rm Bt}\nolimits}
\def\Emb{\mathop{\rm Emb}\nolimits}
\def\Scatt{\mathop{\rm Scatt}\nolimits}
\def\Fin{\mathop{\rm Fin}\nolimits}
\def\sq{\mathop{\rm sq}\nolimits}
\def\asq{\mathop{\rm asq}\nolimits}
\def\ro{\mathop{\rm ro}\nolimits}
\begin{document}
\thispagestyle{plain}
\begin{center}
           {\large \bf \uppercase{Copies of the random graph: the 2-localization}}
\end{center}
\begin{center}
{\bf Milo\v s S.\ Kurili\'c\footnote{Department of Mathematics and Informatics, Faculty of Science, University of Novi Sad,
              Trg Dositeja Obradovi\'ca 4, 21000 Novi Sad, Serbia.
              email: milos@dmi.uns.ac.rs}
and Stevo Todor\v cevi\'c\footnote{Institut de Math\'ematique de Jussieu
              (UMR 7586) Case 247, 4 Place Jussieu, 75252 Paris Cedex, France
              and
              Department of Mathematics, University of Toronto,
              Toronto, Canada M5S 2E4.
              email: stevo.todorcevic@imf-prg.fr and stevo@math.toronto.edu}}
\end{center}
\begin{abstract}
\noindent
Let $G$ be a countable graph containing a copy of the countable random graph (Erd\H{o}s-R\'enyi graph, Rado graph),
$\Emb (G)$ the monoid of its self-embeddings, $\P (G)=\{ f[G]: f\in \Emb (G)\}$ the set of copies
of $G$ contained in $G$, and ${\mathcal I}_G$ the ideal of subsets of $G$ which do not contain a copy of $G$.
We show that the poset $\la \P (G ), \subset\ra$, the algebra $P (G)/{\mathcal I _G}$, and the inverse of the right Green's pre-order $\la \Emb (G),\preceq ^R \ra$
have the 2-localization property.
The Boolean completions of these pre-orders are isomorphic and satisfy the following law:
for each double sequence $[b_{n m }: \langle n , m \rangle \in \o \times \o ]$ of elements of ${\mathbb B}$
$$
\textstyle
\bigwedge _{n \in \o }\;
\bigvee _{m \in \o }\;
b_{n m}
=
\bigvee _{\T \,\in \,\Bt ({}^{<\o }\o)}\;
\bigwedge _{n \in \o }\;
\bigvee _{\f \,\in \,\T \cap {}^{n+1}\o }\;
\bigwedge _{k\leq n}\;
b_{k\f (k)} ,
$$
where $\Bt ({}^{<\o }\o)$ denotes the set of all binary subtrees of the tree ${}^{<\o }\o$.
\vspace{1mm} 
\end{abstract}
\section{Introduction}\label{S1}
The structure considered in this paper is the self-embedding monoid $\Emb (R)$ of the countable random graph
(the Erd\H{o}s-R\'enyi graph \cite{Erdos2}, the Rado graph \cite{Rado}).

More generally, as a part of the investigation of the class of self-embedding monoids of first order structures, in an attempt to classify the pre-orders of the form $\la \Emb (\X), \preceq ^R \ra$, where $\X$ is a structure and
$\preceq ^R$ the right Green's pre-order on the set $\Emb (\X )$ of its self-embeddings (defined by $f\preceq ^R g$ iff $f\circ h=g$, for some $h\in \Emb (\X)$),
one can define two such pre-orders to be equivalent iff the antisymmetric quotients of their inverses  have isomorphic Boolean completions.
It is easy to see \cite{Kmon}
that the antisymmetric quotient of the pre-order $\la \Emb (\X), (\preceq ^R )^{-1} \ra$ is isomorphic to the poset
$\la \P (\X ), \subset \ra$, where $\P (\X )=\{ f[X]: f\in \Emb (\X )\}$ is the set of copies of $\X$ (that is, the set of
domains of the substructures of $\X$ which are isomorphic to $\X$) and that
the isomorphism of the Boolean completions of such posets is the same as their forcing equivalence \cite{Kdif}. Thus, the intended classification is in fact
the classification of the posets of the form $\la \P (\X ), \subset \ra$ determined by their forcing-related properties.
Clearly, this classification induces a coarse classification of structures as well (see \cite{Ktow,Kur1,Kurord,Kurstr,Kuremb}, for countable relational structures).

Concerning the countable ultrahomogeneous relational structures first we mention
the following result from \cite{KurTod} related to the poset of copies of the rational line, $\Q$,
and the corresponding quotient $P(\Q )/\Scatt$, where $\Scatt$ denotes the ideal of scattered suborders of $\Q$:
if $\S$ denotes the Sacks perfect set forcing and sh$(\S )$
the size of the continuum in the Sacks extension, then for
each countable non-scattered linear order $L$ and, in particular, for the rational line, the poset $\P (L)$ is forcing equivalent to
the two-step iteration
$$
\S \ast \pi
$$
where $1_\S \Vdash `` \pi $ is a $\sigma$-closed forcing".
If the equality sh$(\S )=\aleph _1$ (implied by CH)
or PFA holds in the ground model, then the second iterand is forcing equivalent to the poset $(P(\o )/\Fin )^+$ of the
Sacks extension.
Consequently,
$$
\ro \sq \asq \la \Emb (\Q ), (\preceq ^R )^{-1}\ra \cong \ro \sq \P (\Q ) \cong \ro ((P(\Q )/\Scatt )^+)\cong \ro (\S \ast \pi ).
$$
The following similar statement for countable non-scattered graphs (that is, the graphs containing a copy of the Rado graph) was obtained in \cite{KurTodR1}.
\begin{te}  \label{T4}
For each countable non-scattered graph $\la G , \sim \ra$ and, in particular, for the Rado graph, the poset $\P (G)$ is forcing equivalent to
the two-step iteration
\begin{equation}\label{EQ2743}
\P \ast \pi
\end{equation}
where $1_\P \Vdash ``\pi \mbox{ is an }\o\mbox{-distributive forcing}"$ and the poset $\P$ adds a generic real,
has the $\aleph _0$-covering property (thus preserves $\o _1$), has the Sacks property and does not produce splitting reals.
In addition,
\begin{equation}\label{EQ2744}
\ro \sq \asq \la \Emb (G ), (\preceq ^R )^{-1}\ra \cong \ro \sq \P (G) \cong \ro(P (R)/{\mathcal I _R})^+ \cong \ro (\P \ast \pi )
\end{equation}
and these complete Boolean algebras are weakly distributive\footnote{A complete Boolean algebra $\B$ is called {\it weakly distributive}
(or {\it $(\o , \cdot , <\!\!\o)$-distributive}, see \cite{Jech}) iff for each cardinal $\k$ and each matrix
$[b_{n \a }: \langle n , \a \rangle \in \o \times \k ]$ of elements of ${\mathbb B}$ we have
$$
\textstyle
\bigwedge _{n \in \o }\;\;
\bigvee _{\a \in \k }\;\;
b_{n \a}
=
\bigvee _{s: \o \rightarrow [\k ]^{<\o }}\;\;
\bigwedge _{n \in \o }\;\;
\bigvee _{\a \in s(n ) }
b_{n \a} .
$$
}.
\end{te}
We note that the Sacks forcing has all the properties listed in Theorem \ref{T4}:
adds a generic real, has the $\aleph _0$-covering and the Sacks property and does not produce splitting reals.
In the present paper we show that the poset of copies of the Rado graph and, hence, the forcing $\P$ from (\ref{EQ2743})
shares one more property with the Sacks forcing - the 2-localization property. In order to define it
we recall that a subtree $\T$ of the tree ${}^{<\o }\o =\bigcup _{n\in \o}{}^{n}\o$ is called {\it binary} iff each $\f \in \T$ has at most two immediate successors in $\T$; by $\Bt ({}^{<\o }\o)$ we will denote the set of all binary subtrees of ${}^{<\o }\o$. A pre-order $\P$ is said to have
the {\it 2-localization property} iff in each generic extension of the ground model $V$ by $\P$ we have: for each function $x:\o \rightarrow \o$ there is a binary subtree $\T$ of
${}^{<\o }\o$ belonging to $V$ and such that the set of finite approximations of $x$, $\{ x\upharpoonright n : n\in \o\}$, is a branch in $\T$, that is
\begin{equation}\label{EQ2700}\textstyle
1_\P \Vdash _\P \forall x:\check{\o} \rightarrow \check{\o } \;\;
                \exists \T \in  ((\Bt ({}^{<\o }\o))^V )\check{\;}\;\;
                \forall n\in \check{\o } \;\; x\upharpoonright n\in \T .
\end{equation}
It is clear that the 2-localization property implies the Sacks property, but the converse is not true (see \cite{NewRos}).
For a complete Boolean algebra $\B$, the 2-localization property has the following forcing-free translation, similar to the distributivity
laws: for each double sequence $[b_{n m }: \langle n , m \rangle \in \o \times \o ]$ of elements of ${\mathbb B}$
\begin{equation}\label{EQ2745}
\textstyle
\bigwedge _{n \in \o }\;
\bigvee _{m \in \o }\;
b_{n m}
=
\bigvee _{\T \,\in\, \Bt ({}^{<\o }\o)}\;
\bigwedge _{n \in \o }\;
\bigvee _{\f \,\in\, \T \cap {}^{n+1}\o }\;
\bigwedge _{k\leq n}\;
b_{k\f (k)}
\end{equation}
(see Appendix). Thus, the following statement is the main result of the paper.
\begin{te}  \label{T2663}
For each countable non-scattered graph $\la G , \sim \ra$ and, in particular, for the Rado graph, the poset $\P (G)$, and, hence, the
first iterand $\P$ from (\ref{EQ2743}), has the 2-localization property. Consequently, the algebras from (\ref{EQ2744}) satisfy (\ref{EQ2745}).
\end{te}
\section{Preliminaries}\label{S2}
First we introduce a convenient notation. If $\la G, \sim \ra$ is a graph (namely, if $\sim$ is a symmetric and irreflexive binary relation on the set $G$) and $K\subset H \in [G]^{<\o }$, let
$$
G^H_K \!:=\Big\{ v\in G \setminus H : \forall k\in K \, (v\sim k) \; \land \; \forall h\in H\setminus K \, (v\not\sim h)\Big\}.
$$
(Clearly, $G^\emptyset _\emptyset =G$.)
The object of our study is the Rado graph,
characterized as the unique (up to isomorphism) countable graph $\langle R, \sim \rangle $ such that
\begin{equation}\label{EQ2740}
R^H_K\neq \emptyset, \mbox{ whenever } K\subset H\in [R]^{<\o }.
\end{equation}
More information about the Rado graph and the related structures can be found in the survey article \cite{Camer}.
Now we recall some definitions and facts which will be used in the sequel.

\begin{fac}     \label{T2600}
Let $\langle R,\sim \rangle$ be the Rado graph and $\P (R)$ the set of its copies. Then
\begin{itemize}\itemsep=-0.5mm 
\item[\rm (a)] If $F$ is a finite subset of $R$, then $R\setminus F \in \P (R)$;
\item[\rm (b)] If $\{ X_1,  \dots , X_k\} $ is a partition of $R$, then $X_i\in \P (R)$, for some $i\leq k$ (namely, the Rado graph is a strongly indivisible structure);
\item[\rm (c)] If $H$ is a finite subset of $R$, then $\{ H \} \cup \{ R^H_K : K\subset H\}$ is a partition of $R$ and $R^H_K\in \P (R)$, for each $K\subset H$.
\end{itemize}
\end{fac}
\begin{fac}  [\cite{KurTodR1}]   \label{T2620}
Let $H_1$ and $H_2$ be finite subsets of $R$,  $K_1\subset H_1$ and $K_2\subset H_2$.
\begin{itemize}\itemsep=-1mm 
\item[\rm (a)] $R^{H_1}_{K_1}\cap R^{H_2}_{K_2}\neq \emptyset $ if and only if $H_1 \cap K_2 = H_2 \cap K_1$;
\item[\rm (b)] $R^{H_1}_{K_1}\cap R^{H_2}_{K_2}\neq \emptyset$ implies that $R^{H_1}_{K_1}\cap R^{H_2}_{K_2}=R^{H_1\cup H_2}_{K_1\cup K_2}$;
\item[\rm (c)] $R^{H_1}_{K_1}= R^{H_2}_{K_2}$ if and only if $H_1 = H_2$ and $K_1= K_2$;
\item[\rm (d)] $R^{H_1}_{K_1}\subset R^{H_2}_{K_2}$ if and only if $H_1 \supset H_2$,  $K_1\supset K_2$ and $H_2\cap K_1 =K_2$.
\end{itemize}
\end{fac}
If $\la R , \sim \ra$ is the Rado graph and $L\in \P (R)$, we will say that a pair $\CL =\la \Pi , q \ra$ is a {\it labeling} of $L$ iff
\begin{itemize}\itemsep=-1mm \itemindent=2mm
\item[(L1)] $\Pi =\{ L_n : n\in \o\}$ is a partition of the set $L$,
\item[(L2)] $q: \bigcup _{n\in \o }\{ n \} \times P(\bigcup _{i<n} L_i) \rightarrow L$ is a bijection,
\item[(L3)] $L_n =\{ q (n,K) : K\subset \bigcup _{i<n}L_i \}$, for each $n\in \o$,\\[-5mm]
\item[(L4)] $q (n,K) \in  L ^{\bigcup _{i<n}L_i}_K$, for each $n\in \o$ and each $K\subset \bigcup _{i<n}L_i$.
\end{itemize}
Then, clearly, $L_0=\{q(0 ,\emptyset ) \}$,  $|L_0|=1$ and the sets $L_n$ are finite. More precisely, by (L3) we have $|L_n |=m_n$,
where the  integers $m_n$, $n\in \o$, are defined by: $m_0=1$ and $m_n =2^{\sum _{i<n}m_i}$, for $n>0$.
Thus $\la |L_n|:n\in \o \ra =\la 1,2,8, 2^{11}, \dots \ra$.
We note that, by \cite{KurTodR1}, each copy $L\in \P (R)$ has infinitely many labelings.
For convenience, instead of $q(n,K)$ we will write $q ^{\bigcup _{i<n}L_i}_K$ and the labeling $\CL$
will be denoted by
$$
\textstyle \Big\la \{ L_n : n\in \o \}, \{ q ^{\bigcup _{i<n}L_i}_K : n\in \o \land K\subset \bigcup _{i<n}L_i\}\Big\ra.
$$
If $\D= \la \D _n :n\in \o \ra$ is a sequence of subsets of $\P (R)$, then  a copy $L\in \P (R)$ will be called
a {\it fusion} of the sequence $\D$ if and only if there exists a labeling
$\la \{ L_n : n\in \o \}, \{ q ^{\bigcup _{i<n}L_i}_K : n\in \o \land K\subset \bigcup _{i<n}L_i\}\ra$
of $L$ such that
\begin{equation}\label{EQ2602}\textstyle
\forall n\in \o \;\; \forall  K \subset \bigcup _{i<n}L_i\;\;  \exists D \in \D _n \;\;L^{\bigcup _{i<n}L_i}_K \subset D .
\end{equation}
\begin{fac} [\cite{KurTodR1}]  \label{T2659}
If  $\D= \la \D _n \!:n\in \o \ra$ is a sequence subsets of $\P (R)$ which are dense below $A\in \P (R)$,
then the set ${\mathcal F}=\{L : L \mbox{ is a fusion of }\D \}$ is dense below $A$.
\end{fac}
Using a labeling
$\CL =\la \{ L_n : n\in \o \}, \{ q ^{\bigcup _{i<n}L_i}_K : n\in \o \land K\subset \bigcup _{i<n}L_i\}\ra$ of a copy $L\in \P (R)$
we define a partial ordering $\leq _{L,\CL}$  on $L$ by:
\begin{eqnarray}\label{EQ2622}\textstyle
q ^{\bigcup _{i<m}L_i}_{K '} \leq _{L,\CL} \; q ^{\bigcup _{i<n}L_i}_{K''} & \Leftrightarrow & L^{\bigcup _{i<m}L_i}_{K '} \subset  L ^{\bigcup _{i<n}L_i}_{K''}\\
                                                                   & \Leftrightarrow & \label{EQ2605}\textstyle m\geq n \;\;\land \;\; K' \cap \bigcup _{i<n}L_i =K''
\end{eqnarray}
(see Fact \ref{T2620}(c) and (d)). Writing shortly $\leq _L$ instead of $\leq _{L,\CL}$ we have
\begin{fac} [\cite{KurTodR1}]  \label{T2627}
$\la L, \leq _L  \ra$ is a finitely branching reversed tree without minimal nodes, with the top $q^\emptyset _\emptyset $,
and the set $L_n$ is its $n$-th level. For each $n\in\o$ and $K \subset \bigcup _{i<n}L_i$ we have
$$
\Big( -\infty , q ^{\bigcup _{i<n}L_i}_K \Big]_{\la L, \leq _L  \ra}= L^{\bigcup _{i<n}L_i}_K .
$$
\end{fac}
\section{The 2-localization property}\label{S3}
In this section we prove Theorem \ref{T2663}.
In fact, if $\la G , \sim \ra$ is a countable graph containing a copy of the Rado graph, then these two structures are equimorphic and, by \cite{Kdif},
the corresponding posets of copies $\P (G)$ and $\P (R)$ are forcing equivalent. So it is sufficient to prove the theorem assuming that $\la G , \sim \ra$ is the Rado graph and we will do this in Theorem \ref{T2638}. We start with two auxiliary statements.
\begin{lem}  \label{T2637}
Let $k_n\in \o $, $n\in \o$, where $k_0=0$ and $k_{n+1}=k_n +2^{n+1}-1$, and let
${}^n 2= \{ \f ^n_i : i<2^n \}$ be an enumeration of the $n$-th level of the binary tree ${}^{<\o }2=\bigcup _{n\in \o}{}^n 2$
induced by the lexicographic order (for example: $\f ^2 _0= 00$, $\f ^2 _1= 01$, $\f ^2 _2= 10$, and $\f ^2 _3= 11$). Then
$$
\Delta = \Big\{ \f ^{n+1}_i {}^{\smallfrown} \underbrace{00\dots 0}_{k_n + \,i} : n\in \o \;\land \; i<2^{n+1}\Big\}
$$
is a dense subset of the reversed tree $\la{}^{<\o }2, \supset \ra$ and $|\Delta \cap {}^n 2|=1$, for each $n\in \N$.
\end{lem}
\dok
The density of $\Delta$ is evident. Since the lengths of the elements of $\Delta$ below the $(n+1)$-th level of ${}^{<\o }2$ run from $n+1 +k_n$ to
$n+1 +k_n +2^{n+1}-1=n+1 +k_{n+1}$, while the lengths of the elements of $\Delta$ below the $(n+2)$-th level of ${}^{<\o }2$ start from $n+2+k_{n+1}$, the second statement is true as well.
\hfill $\Box$
\begin{lem}  \label{T2639}
If $g:\la {}^{<\o }2 ,\supset \ra \rightarrow \la {}^{<\o }\o ,\supset \ra$ is an embedding (that is, an injective strong homomorphism), then the set
$g[{}^{<\o }2]\uparrow =\{ \p \in {}^{<\o }\o :\exists \f \in {}^{<\o }2\; g(\f )\supset \p \}$
is a binary reversed subtree of the reversed tree $\la {}^{<\o }\o ,\supset \ra$.
\end{lem}
\dok
Suppose that $\p \in g[{}^{<\o }2]\uparrow$ has 3 immediate predecessors $\p _i $, $i<3$, in  $g[{}^{<\o }2]\uparrow$.
Let $\f _i$, $i<3$, be the elements of ${}^{<\o }2$ of minimal length satisfying $g(\f _i)\supset \p _i$.
Then $g(\f _0 \cap \f _1)\subset \p _0 \cap \p _1=\p$.
Let $\f $ be the element of ${}^{<\o }2$ of maximal length satisfying $\p \supset g(\f )$.
Suppose that for some $i<3$ there is $\e\in {}^{<\o }2$ such that $\f _i \supsetneq \e \supsetneq \f$.
Then $g(\f _i )\supsetneq g(\e )\supsetneq g(\f )$ and, since $g(\f _i )\supsetneq g(\e ) , \p$ and $\la {}^{<\o }\o ,\supset \ra$ is a reversed tree,
$g(\e )$ and $\p$ are comparable. But $\p \supset g(\e )$ is not true, by the maximality of $\f$ and $g(\e )\supsetneq \p$ is not true, by the minimality of $\f _i$.
Thus $\f _i$, $i<3$, are immediate predecessors of $\f$ in ${}^{<\o }2$, which is impossible.
\hfill $\Box$
\begin{te}  \label{T2638}
The partial ordering $\P (R)$ has the 2-localization property.
\end{te}
\dok
Let $G$ be a $\P (R)$-generic filter over the ground model $V$ and $x\in V_{\P (R)}[G]$, where $x : \o \rightarrow \o$.
Let $\t$ be a $\P (R)$-name such that $x=\t _G$ and let $A\in G$, where $A\Vdash \t : \check{\o }\rightarrow \check{\o }$.
First we define the sets $\D_0, \D _1\subset \P (R)$ by
$$
\D _0 =\{ B\in \P (A) : \forall n\in \o \;\; \exists m\in \o \;\; B\Vdash \t (\check{n})= \check{m}\} \mbox{ and }
$$
$$
\D _1=\{ C\in \P (A) : \forall B\in \P (C)\;\; \exists n\in \o \;\; \forall m\in \o \;\; \neg B\Vdash \t (\check{n})= \check{m}\}
$$
and show that, in the poset $\la \P (R), \subset \ra$, the union $\D =\D _0 \cup \D_1$ is dense below $A$.
So, let $D\in \P (R)$ and $D\subset A$. If $D\in \D _1$, we are done; otherwise, there is $B\in \P (D)$ such that $B\in \D _0$, thus $\D \ni B\subset D$ and
$\D $ is dense below $A$ indeed.

Since $A\in G$, the genericity of $G$ implies $G\cap \D\neq \emptyset$ and we have two cases.

\vspace{2mm}

\noindent
{\bf Case 1:} $G\cap \D _0\neq \emptyset$. Then $B\in G$, for some $B\in \D _0$,
and, by the definability of the forcing relation, the set $f=\{ \la n,m\ra \in \o \times \o:  B\Vdash \t (\check{n})= \check{m}\}$ belongs to $V$.
In order to prove that $B\Vdash \t =\check{f}$ we assume that $H$ is a $\P (R)$-generic filter over $V$ containing $B$ and show that $\t _H=f$.
If $\la n,m \ra \in \t _H$, then, since $B\in \D _0$, there is $p\in \o$ such that $B\Vdash \t (\check{n})= \check{p}$, which implies $\la n,p \ra \in \t _H$;
thus, since $B\Vdash \t : \check{\o }\rightarrow \check{\o }$, we have $p=m$ and, hence $\la n,m \ra \in f$. Conversely, if $\la n,m \ra \in f$, then
$B\Vdash \t (\check{n})= \check{m}$ and, since $B\in H$, we have $\la n,m \ra \in \t _H$. So $B\Vdash \t =\check{f}$ and, by the assumption, $B\in G$,
which implies $\t _G =f \in V$. Thus for $\T =\{ f\upharpoonright n:n\in \o\} $ we have $\T \in \Bt ({}^{<\o}\o )$ and
$\t _G \upharpoonright n \in \T$, for all $n\in \o$.

\vspace{2mm}

\noindent
{\bf Case 2:} $G\cap \D _0= \emptyset$. Then, by the density of $\D$, there is $C\in G\cap \D _1$ so we have
\begin{equation}\label{EQ2701}\textstyle
\forall B\in \P (C) \;\; \exists n\in \o \;\; \forall m\in \o \;\; \neg B\Vdash \t (\check{n})= \check{m}
\end{equation}
and, since $C\in G$, it is sufficient to prove that
\begin{equation}\label{EQ2702}\textstyle
\forall B\in \P (C) \;\; \exists D\in \P (B) \;\; \exists \T \in \Bt ({}^{<\o}\o ) \;\;  \forall n\in \o \;\;D\Vdash \t \upharpoonright \check{n}\in  \check{\T}.
\end{equation}
Let $B\in \P (C)$. Then, since $B\subset C\subset A$, we have $B\Vdash \t : \check{\o }\rightarrow \check{\o }$ and, hence,
$B\Vdash  \forall n\in \check{\o } \; \exists m \in \check{\o }\; \t (n)=m$. Thus, for each $n\in \o$ we have:
for each $B'\in \P (B)$ there are $D\in \P (B')$ and $m \in \o$ such that $D\Vdash \t (\check{n})=\check{m}$. This means that the sets
$\D _n :=\{ D\in \P (B): \exists m \in \o \;\; D\Vdash \t (\check{n})=\check{m}\}$, $n\in \o$, are dense below $B$.
By Fact \ref{T2659}, the set ${\mathcal F}$ of fusions of the sequence $\la \D _n : n\in \o\ra$  is dense below $B$
and, hence, there is a fusion $L \in {\mathcal F}$ such that $L\in \P (B)$. So there is a labeling
$\la \{ L_n : n\in \o \}, \{ q ^{\bigcup _{i<n}L_i}_K : n\in \o \land K\subset \bigcup _{i<n}L_i\}\ra$ of $L$ such that, by (\ref{EQ2602}), for
$n\in \o$ and $K \subset \bigcup _{i<n}L_i$ there is $D \in \D _n $ such that $L^{\bigcup _{i<n}L_i}_K \subset D $ and, hence, there is (clearly unique)
$m \in \o$ satisfying $L^{\bigcup _{i<n}L_i}_K \Vdash \t (\check{n})=\check{m}$.
Thus we obtain an indexed family of integers $\{ m^n_K :n\in \o \land K\subset \bigcup _{i<n}L_i\}$ such that
\begin{equation}\label{EQ2703}\textstyle
\forall n\in \o \;\; \forall K\subset \bigcup _{i<n}L_i\;\;  L^{\bigcup _{i<n}L_i}_K\Vdash \t (\check{n})=\check{m^n_K}.
\end{equation}
Let $\leq _L$ be the reversed tree order on $L$ defined by (\ref{EQ2622}) or, equivalently, by (\ref{EQ2605}).

\begin{cla}\label{T2651}
For each $n\in \o$ and each $K\subset \bigcup _{i<n}L_i$ there are $n_1>n$ and $K',K''\subset \bigcup _{i<n_1}L_i$ such that
\begin{equation}\label{EQ2704}\textstyle
q^{\bigcup _{i<n_1}L_i}_{K'} , q^{\bigcup _{i<n_1}L_i}_{K''}<_L q^{\bigcup _{i<n}L_i}_K \;\;\land\;\;
m^{n_1}_{K'}\neq m^{n_1}_{K''} .
\end{equation}
\end{cla}
\dok
Suppose, contrary to our claim, that there are $n_0\in \o$ and $K_0\subset \bigcup _{i<n_0}L_i$ such that for  each $n>n_0$
and each $K',K''\subset \bigcup _{i<n}L_i$ we have (see (\ref{EQ2605}))
\begin{equation}\label{EQ2705}\textstyle
K' \cap \bigcup _{i<n_0}L_i =K'' \cap \bigcup _{i<n_0}L_i =K_0 \;\;\Rightarrow \;\; m^n_{K'}= m^n_{K''}.
\end{equation}
Then for  each $n>n_0$ there is $m(n)\in \o$ such that for each $K\subset \bigcup _{n_0\leq i<n}L_i$ we have $m^n_{K_0\cup K}=m(n)$
which, by (\ref{EQ2703}), implies
$L^{\bigcup _{i<n}L_i}_{K_0\cup K}=(L^{\bigcup _{i<n_0}L_i}_{K_0} )^{\bigcup _{n_0\leq i <n}L_i}_K\Vdash \t (\check{n})= \check{m(n)}$.
By Fact \ref{T2600}(c) $\{ (L^{\bigcup _{i<n_0}L_i}_{K_0} )^{\bigcup _{n_0\leq i <n}L_i}_K :K\subset \bigcup _{n_0\leq i<n}L_i\}$ is an antichain in $\P (R)$,
maximal below $L^{\bigcup _{i<n_0}L_i}_{K_0}$, thus $L^{\bigcup _{i<n_0}L_i}_{K_0} \Vdash \t (\check{n})= \check{m(n)}$ and this holds for all $n>n_0$.
For $n\leq n_0$, let $m(n)=m^{\bigcup _{i<n}L_i}_{K_0 \cap \bigcup _{i<n}L_i}$; then,
since $L^{\bigcup _{i<n_0}L_i}_{K_0}\subset L^{\bigcup _{i<n}L_i}_{K_0 \cap \bigcup _{i<n}L_i}$, by (\ref{EQ2703})
we have $L^{\bigcup _{i<n_0}L_i}_{K_0} \Vdash \t (\check{n})= \check{m(n)}$. Thus
$L^{\bigcup _{i<n_0}L_i}_{K_0}\in \P (C)$ and $L^{\bigcup _{i<n_0}L_i}_{K_0} \Vdash \t (\check{n})= \check{m(n)}$, for each $n\in \o$, which is not true by (\ref{EQ2701}).
\kdok

\noindent
{\bf Strategy.}
Roughly speaking, in order to prove (\ref{EQ2702}), using the fusion $L\subset B$ obtained above,
its labeling and the corresponding order $\leq _L$, we will construct:

- a copy $\L $ of the reversed binary tree $\la {}^{<\o }2, \supset \ra$ in the reversed tree $\la L , \leq _L \ra$,

- a copy $D$ of the Rado graph contained in $\L$,

- an embedding $F$ of $\L$ in the reversed tree $\la {}^{<\o }\o, \supset \ra$ such that for the upwards closure $\T $ of $F[\L ]$
we have $D\Vdash \{ \t \upharpoonright n : n\in \check{\o }\}\subset \check{\T }$.


\vspace{2mm}

\noindent
{\bf Construction.}
Let $\Delta=\{ \p _n :n\in \N \}$ be the enumeration of the set $\Delta\subset {}^{<\o}2$, defined in Lemma \ref{T2637}, such that $\Delta\cap {}^n 2 =\{ \p _n \}$.
For convenience let $\L _0 := \{ q ^\emptyset _\emptyset\}$,  $l_0 :=0$ and $K_\emptyset :=\emptyset$.
Using recursion, for each $n\in \N$ we define $\L _n , f_n , l_n , d_n , \F _n$, $\{ n_\f :\f \in {}^{n-1}2\}$ and $\{ K_\f :\f \in {}^{n}2\}$ such that for each
$n\in \N$ we have:
\begin{itemize}\itemsep=-1mm \itemindent=5mm
\item[(i)]  $K_{\f ^\smallfrown j}\subset \bigcup _{i<n_\f}L_i$, for each $\f \in {}^{n-1}2$ and $j\in 2$;
\item[(ii)] $n_\f \leq n_{\p _n \upharpoonright (n-1)}< n_{\p _n \upharpoonright (n-1)}+1 =l_n <n_{\f ^\smallfrown j}$, for $\f \in {}^{n-1}2$ and $j\in 2$;
\item[(iii)] $q^{\bigcup _{i<n_{\f ^\smallfrown 0}}L_i}_{K_{\f ^{\smallfrown} 0 {}^\smallfrown j}}
             <_L   q^{\bigcup _{i<l_{n} }L_i}_{K_{\f ^{\smallfrown} 0 }}$, for each $\f \in {}^{n-1}2$ and $j\in 2$;
\item[(iv)] $q^{\bigcup _{i<n_{\f ^\smallfrown 1}}L_i}_{K_{\f ^{\smallfrown} 1 {}^\smallfrown j}}
            <_L   q^{\bigcup _{i<l_n} L_i}_{K_{\f ^{\smallfrown} 1 } \cup \{ d_n \}}$, for each $\f \in {}^{n-1}2$ and $j\in 2$;
\item[(v)]  $m^{n_\f }_{K_{\f ^{\smallfrown} 0}}\neq m^{n_\f }_{K_{\f ^{\smallfrown} 1}}$, for each $\f \in {}^{n-1}2$;
\item[(vi)] $\L _n =\{ q ^{\bigcup _{i<n_\f}L _i}_{K_{\f ^{\smallfrown}j}} : \f\in {}^{n-1}2 \land j\in 2   \} $;
\item[(vii)] $f_n :\la \bigcup _{i\leq n}{}^i 2 , \supset \ra \rightarrow \la \bigcup _{i\leq n} \L _i , \leq _L\ra$ is an isomorphism, where
             $$f_n(\emptyset)= q ^\emptyset _\emptyset \;\;\mbox{ and } \;\;f_n(\f ^{\smallfrown}j)=q ^{\bigcup _{i<n_\f}L _i}_{K_{\f ^{\smallfrown}j}};$$
\item[(viii)] $f_m\subset f_n$, for all $m\leq n$;
\item[(ix)] $d_n=f_n (\p _n)$;
\item[(x)]  $\F _n =\{q^{\bigcup _{i<l_{n} }L_i}_{K_{\f ^{\smallfrown} 0 }}: \f \in {}^{n-1}2\}
            \cup  \{ q^{\bigcup _{i<l_n} L_i}_{K_{\f ^{\smallfrown} 1 } \cup \{ d_n \}}: \f \in {}^{n-1}2\}$.
\end{itemize}
\begin{cla}\label{T2652}
The recursion works.
\end{cla}
\dok
By Claim \ref{T2651} (for $n=0$ and $K=\emptyset$) there are $n_\emptyset >0$ and $K_0,K_1 \subset \bigcup _{i< n_\emptyset}L_i$ such that
$q^{\bigcup _{i<n_\emptyset}L_i}_{K_0} , q^{\bigcup _{i<n_\emptyset }L_i}_{K_1}<_L q^\emptyset_\emptyset$ and
$m^{n_\emptyset}_{K_0}\neq m^{n_\emptyset}_{K_1}$. Let $\L _1 =\{ q^{\bigcup _{i<n_\emptyset}L_i}_{K_j}:j\in 2\}$ let
$f_1 : \la {}^0 2 \cup {}^1 2, \supset \ra \rightarrow \la \L _0 \cup \L_1 ,\leq _L\ra$, where
$f_1(\emptyset)=q^\emptyset _\emptyset $ and $f_1(\la j\ra )=q^{\bigcup _{i<n_\emptyset}L_i}_{K_j}$, for $j\in 2$,
let $l_1=n_\emptyset +1$, $d_1=f_1 (\p _1)=f_1 (\la 0\ra)=q^{\bigcup _{i<n_\emptyset}L_i}_{K_0}$
and $\F _1 =\{ q^{\bigcup _{i<l_1}L_i}_{K_0}, q^{\bigcup _{i<l_1}L_i}_{K_1 \cup \{ d_1 \}}\}$. It is easy to check that conditions
(i) - (x) are satisfied.

Suppose that the objects $\L _i$, $f_i $, $l_i $, $d_i $, and $\F _i$, for $i\leq n$, $\{ n_\f :\f \in {}^{\leq n-1 }2\}$ and $\{ K_\f: \f \in {}^{\leq n }2\}$
satisfy conditions (i) - (x).

First, for $\f \in {}^{n-1 }2$ and $k,j\in 2$, we define $n_{\f ^{\smallfrown} k}$ and $K_{\f ^\smallfrown k{}^\smallfrown j}$. By (ii) and (i) we have
\begin{equation}\label{EQ2708}\textstyle
l_n >n_\f \;\;\mbox{ and } \;\;
K_{\f ^\smallfrown k}\subset \bigcup _{i<n_\f}L_i \subset \bigcup _{i<l_n}L_i
\end{equation}
Since all the sequences from $\Delta$ have 0 at the end we have  $\p _{n+1}=\f _0 ^{\smallfrown}k_0 {}^{\smallfrown}0$, for some $\f _0\in {}^{n-1} 2$ and some $k_0\in 2$, and we distinguish the following four cases.
\begin{itemize}
\item[{\sc (i)}] $\f ^{\smallfrown} k \neq \f _0 ^{\smallfrown}k_0$ and $k=0$. By (\ref{EQ2708}) and  Claim \ref{T2651}
(applied to $l_n$ and $K_{\f ^\smallfrown 0}$) there are
\begin{equation}\label{EQ2706}\textstyle
n_{\f ^{\smallfrown} 0}> l_n \;\;\mbox{ and } \;\; K_{\f ^\smallfrown 0{}^\smallfrown j}\subset \bigcup _{i<n_{\f ^{\smallfrown} 0}}L_i , \mbox{ for } j\in 2,
\mbox{ such that}
\end{equation}
\begin{equation}\label{EQ2707}\textstyle
q^{\bigcup _{i<n_{\f ^{\smallfrown} 0}}L_i}_{K_{\f ^\smallfrown 0{}^\smallfrown j}}
<_L q^{\bigcup _{i<l_n}L_i}_{K_{\f ^\smallfrown 0}}
<_L q^{\bigcup _{i<n_{\f }}L_i}_{K_{\f ^\smallfrown 0}}
\;\;\mbox{ and } \;\;
m^{n_{\f ^{\smallfrown} 0}}_{K_{\f ^\smallfrown 0{}^\smallfrown 0}}\neq m^{n_{\f ^{\smallfrown} 0}}_{K_{\f ^\smallfrown 0 {}^\smallfrown 1}} .
\end{equation}
\item[{\sc (ii)}] $\f ^{\smallfrown} k \neq \f _0 ^{\smallfrown}k_0$ and $k=1$. Then, by (ix), (vii), (vi) and (ii),
\begin{equation}\label{EQ2714}\textstyle
d_n=f_n (\p _n)\in \L _n \subset \bigcup _{\f '\in {}^{n-1} 2}L_{n_{\f '}} \subset \bigcup _{i<l_n}L_i ,
\end{equation}
so, by  (\ref{EQ2708}) we have $K_{\f ^\smallfrown 1} \cup \{ d_n \}\subset \bigcup _{i<l_n}L_i $ and by
Claim \ref{T2651} (applied to $l_n$ and $K_{\f ^\smallfrown 1}\cup \{ d_n \}$) there are
\begin{equation}\label{EQ2709}\textstyle
n_{\f ^{\smallfrown} 1}> l_n \;\;\mbox{ and } \;\; K_{\f ^\smallfrown 1{}^\smallfrown j}\subset \bigcup _{i<n_{\f ^{\smallfrown} 1}}L_i , \mbox{ for } j\in 2,
\mbox{ such that}
\end{equation}
\begin{equation}\label{EQ2710}\textstyle
q^{\bigcup _{i<n_{\f ^{\smallfrown} 1}}L_i}_{K_{\f ^\smallfrown 1{}^\smallfrown j}}
<_L q^{\bigcup _{i<l_n}L_i}_{K_{\f ^\smallfrown 1}\cup \{ d_n \} }
\;\;\mbox{ and } \;\;
m^{n_{\f ^{\smallfrown} 1}}_{K_{\f ^\smallfrown 1{}^\smallfrown 0}}\neq m^{n_{\f ^{\smallfrown} 1}}_{K_{\f ^\smallfrown 1 {}^\smallfrown 1}} .
\end{equation}
\end{itemize}
So $n_{\f ^{\smallfrown} k}$ and $K_{\f ^\smallfrown k{}^\smallfrown j}$ are defined for $\f ^{\smallfrown} k \neq \f _0 ^{\smallfrown}k_0$
and now we define $n_{\f _0 ^{\smallfrown}k_0}$ and $K_{\f _0^\smallfrown k_0{}^\smallfrown j}$, for $j\in 2$. Let
$n^*=\max \{ n_{\f ^\smallfrown k}  : \f \in {}^{n-1} 2 \land k\in 2 \land \f  ^{\smallfrown} k \neq \f _0 ^{\smallfrown}k_0 \} $.
By (\ref{EQ2706}) and (\ref{EQ2709}) we have $n^*>l_n$
thus, by (\ref{EQ2708}),
\begin{equation}\label{EQ2713}\textstyle
K_{\f _0 ^\smallfrown k_0}\subset \bigcup _{i<l_n}L_i \subset \bigcup _{i<n^*}L_i .
\end{equation}
\begin{itemize}
\item[{\sc (iii)}] $\f ^{\smallfrown} k = \f _0 ^{\smallfrown}k_0 $ and $k=0$.  By (\ref{EQ2713}) and
Claim \ref{T2651} (applied to $n^*$ and $K_{\f _0 ^\smallfrown 0}$) there are
\begin{equation}\label{EQ2711}\textstyle
n_{\f _0 ^{\smallfrown} 0}> n^* > l_n >n_{\f _0}
\;\;\mbox{ and } \;\;
K_{\f _0 ^\smallfrown 0{}^\smallfrown j}\subset \bigcup _{i<n_{\f _0 ^{\smallfrown} 0}}L_i , \mbox{ for } j\in 2,
\mbox{ where}
\end{equation}
\begin{equation}\label{EQ2712}\textstyle
q^{\bigcup _{i<n_{\f _0 ^{\smallfrown} 0}}L_i}_{K_{\f _0 ^\smallfrown 0{}^\smallfrown j}}
<_L q^{\bigcup _{i<l_n}L_i}_{K_{\f _0 ^\smallfrown 0}}
<_L q^{\bigcup _{i<n_{\f _0}}L_i}_{K_{\f _0^\smallfrown 0}}
\;\;\mbox{ and } \;\;
m^{n_{\f _0 ^{\smallfrown} 0}}_{K_{\f _0 ^\smallfrown 0{}^\smallfrown 0}}\neq m^{n_{\f _0^{\smallfrown} 0}}_{K_{\f _0^\smallfrown 0 {}^\smallfrown 1}} .
\end{equation}

\item[{\sc (iv)}] $\f ^{\smallfrown} k = \f _0 ^{\smallfrown}k_0 $ and $k=1$.
For $\f ' \in {}^{n-1}2$ by (ii) we have $n_{\f '} < l_n <n^*$ and, by (\ref{EQ2714})and (\ref{EQ2713}),
$K_{\f _0 ^{\smallfrown} 1}\cup \{ d_n \}\subset \bigcup _{i<n^*}L_i$.
So, by Claim \ref{T2651} (applied to $n^*$ and $K_{\f _0 ^\smallfrown 1}\cup \{ d_n \}$) there are
\begin{equation}\label{EQ2715}\textstyle
n_{\f _0 ^{\smallfrown} 1}> n^* > l_n
\;\;\mbox{ and } \;\;
K_{\f _0 ^\smallfrown 1{}^\smallfrown j}\subset \bigcup _{i<n_{\f _0 ^{\smallfrown} 1}}L_i , \mbox{ for } j\in 2,
\mbox{ such that}
\end{equation}
\begin{equation}\label{EQ2716}\textstyle
q^{\bigcup _{i<n_{\f _0 ^{\smallfrown} 1}}L_i}_{K_{\f _0 ^\smallfrown 1{}^\smallfrown j}}
<_L q^{\bigcup _{i<n^*}L_i}_{K_{\f _0 ^\smallfrown 1}\cup \{ d_n \}}
<_L q^{\bigcup _{i<l_n}L_i}_{K_{\f _0^\smallfrown 1}\cup \{ d_n \}}
\;\;\mbox{ and } \;\;
m^{n_{\f _0 ^{\smallfrown} 1}}_{K_{\f _0 ^\smallfrown 1{}^\smallfrown 0}}\neq m^{n_{\f _0^{\smallfrown} 1}}_{K_{\f _0^\smallfrown 1 {}^\smallfrown 1}} .
\end{equation}
\end{itemize}
So we have defined $n_{\f ^{\smallfrown} k}$ and $K_{\f  ^\smallfrown k{}^\smallfrown j}$, for $\f \in {}^{n-1 }2$ and $k,j\in 2$. Now we define
$$\textstyle
\L _{n+1}:=\Big\{ q^{\bigcup _{i<n_{\f ^{\smallfrown} k}}L_i}_{K_{\f  ^\smallfrown k{}^\smallfrown j}} : \f \in {}^{n-1 }2 \land k,j \in 2 \Big\}
$$
and (vi) is true.
(i) follows from (\ref{EQ2706}), (\ref{EQ2709}), (\ref{EQ2711}) and (\ref{EQ2715}) and
(iii), (iv) and (v) follow from (\ref{EQ2707}), (\ref{EQ2710}), (\ref{EQ2712}) and (\ref{EQ2716}).
Recall that $\p _{n+1}=\f _0 ^{\smallfrown}k_0 ^{\smallfrown}0$ and let us define
$$
l_{n+1}=n_{\f _0 ^{\smallfrown}k_0}+1.
$$
Then,  by (\ref{EQ2711}) and (\ref{EQ2715}), for $\f \in {}^{n-1 }2$ and $k\in 2$ satisfying $\f ^{\smallfrown} k \neq \f _0 ^{\smallfrown}k_0$ we have
\begin{equation}\label{EQ2721}\textstyle
l_{n+1}= n_{\p _{n+1}\upharpoonright n}+1 >n_{\p _{n+1}\upharpoonright n}= n_{\f _0 ^{\smallfrown}k_0}>n^* \geq n_{\f ^{\smallfrown} k}
\end{equation}
and, by (\ref{EQ2706}), (\ref{EQ2709}), (\ref{EQ2711}) and (\ref{EQ2715}), $n_{\f ^{\smallfrown} k}>l_n$, for all $\f \in {}^{n-1 }2$ and $k\in 2$.
Thus (ii) is true as well.

Let $f_{n+1}: \bigcup _{i\leq n+1}{}^i 2\rightarrow \bigcup _{i\leq n+1}\L _i$  be defined by: $f_{n+1}\upharpoonright \bigcup _{i\leq n}{}^i 2=f_n$
(thus (viii) is true) and
\begin{equation}\label{EQ2717}\textstyle
f_{n+1}(\f  ^\smallfrown k{}^\smallfrown j)=
q^{\bigcup _{i<n_{\f ^{\smallfrown} k}}L_i}_{K_{\f  ^\smallfrown k{}^\smallfrown j}}, \mbox{ for }\f \in {}^{n-1 }2 \mbox{ and }k,j \in 2.
\end{equation}
In order to prove that the mapping $f_{n+1}$ is an isomorphism, first we show that
\begin{equation}\label{EQ2718}\textstyle
\forall \f \in {}^{n-1 }2 \;\; \forall k,j \in 2 \;\; \mbox{}
q^{\bigcup _{i<n_{\f ^{\smallfrown} k}}L_i}_{K_{\f  ^\smallfrown k{}^\smallfrown j}} <_L
q^{\bigcup _{i<n_{\f }}L_i}_{K_{\f  ^\smallfrown k}}.
\end{equation}
In cases {\sc (i)} and {\sc (iii)} this follows from (\ref{EQ2707}) and (\ref{EQ2712}).
For case {\sc (ii)}, by (\ref{EQ2710}) and (\ref{EQ2605})
it is sufficient to prove that $(K _{\f ^{\smallfrown} 1}\cup \{ d_n \}) \cap \bigcup _{i<n_{\f }}L_i =K _{\f ^{\smallfrown} 1}$, and, by (i), this
will follow from $d_n \not\in \bigcup _{i<n_{\f }}L_i$. By (ix) and (vii) we have
\begin{equation}\label{EQ2719}\textstyle
d_n
= f_n (\p _n)
= f_n ((\p _n\upharpoonright (n-1))^{\smallfrown}0)
= q^{\bigcup _{i<n_{\p _n\upharpoonright (n-1)}}L_i}_{K_{(\p _n\upharpoonright (n-1))^{\smallfrown}0}}
\in L_{n_{\p _n\upharpoonright (n-1)}}
\end{equation}
but, by (ii), $n_{\p _n\upharpoonright (n-1)}\geq n_\f$ and $d_n \not\in \bigcup _{i<n_{\f }}L_i$ indeed.
For case {\sc (iv)}, by (\ref{EQ2716}) and (\ref{EQ2605})
we have to prove that $(K _{\f _0 ^{\smallfrown} 1}\cup \{ d_n \}) \cap \bigcup _{i<n_{\f _0}}L_i =K _{\f_ 0 ^{\smallfrown} 1}$, and, by (i), this
will follow from
$d_n \not\in \bigcup _{i<n_{\f _0}}L_i$. By (\ref{EQ2719}) we have
$d_n \in L_{n_{\p _n\upharpoonright (n-1)}}$
and, by (ii), $n_{\p _n\upharpoonright (n-1)}\geq n_{\f _0}$. Thus $d_n \not\in \bigcup _{i<n_{\f _0 }}L_i$ and the proof of (\ref{EQ2718}) is finished.

By (vi) and (vii) $\L _n$ is an antichain in $\la L,\leq _L\ra$ and, by (\ref{EQ2707}), (\ref{EQ2710}), (\ref{EQ2712}) and (\ref{EQ2716}),
using Fact \ref{T2627} and (\ref{EQ2703}) we conclude that
$q^{\bigcup _{i<n_{\f ^{\smallfrown} k}}L_i}_{K_{\f  ^\smallfrown k{}^\smallfrown j}} $, $j\in 2$,
are $\leq _L$-incompatible elements below
$q^{\bigcup _{i<n_{\f }}L_i}_{K_{\f  ^\smallfrown k}}$.
So  $f_{n+1}$ is an isomorphism and (vii) is true.

Finally we define $d_{n+1}$ and $\F _{n+1}$ by
\begin{equation}\label{EQ2720}\textstyle
d_{n+1}
= f_{n+1} (\p _{n+1})
= f_{n+1} (\f _0 ^{\smallfrown} k_0{}^{\smallfrown}0 )
= q^{\bigcup _{i<n_{\f _0 ^{\smallfrown} k_0}}L_i}_{K_{\f _0 ^{\smallfrown} k_0{}^{\smallfrown}0}}
\end{equation}
$$\textstyle
\F _{n+1} =\bigcup _{\f \in {}^{n-1}2 \;\land \; k\in 2}
\Big\{
q^{\bigcup _{i<l_{n+1} }L_i}_{K_{\f ^{\smallfrown} k {}^{\smallfrown}0}},
q^{\bigcup _{i<l_{n+1}} L_i}_{K_{\f ^{\smallfrown} k {}^{\smallfrown}1 } \cup \{ d_{n+1} \}}
\Big\}.
$$
Then (ix) is true. By (\ref{EQ2720})
$d_{n+1} \in L_{n_{\f _0 ^{\smallfrown} k_0}}$ and, by (\ref{EQ2721}) $n_{\f _0 ^{\smallfrown} k_0}<l_{n+1}$. Thus $\F _{n+1}$ is well defined
and (x) is true. Claim \ref{T2652} is proved.
\hfill $\Box$
\begin{cla}\label{T2661}
Let $\L :=\bigcup _{n\in \o}\L _n$, $f:=\bigcup _{n\in \N }f _n$ and $D:=\{ d_n : n\in \N \}$. Then
\begin{itemize}\itemsep=-0.5mm 
\item[\rm (a)] The mapping $f: \la {}^{<\o }2 ,\supset \ra \rightarrow \la \L ,\leq _L \ra$ is an isomorphism;
\item[\rm (b)] $D$ is a dense subset of the binary reversed tree $\la \L ,\leq _L \ra$.
\end{itemize}
\end{cla}
\dok
(a) follows from (vii) and (viii).

(b) By (ix) we have $D=f[\Delta ]$ and, by Lemma \ref{T2637}, $\Delta$ is a dense subset of the binary reversed tree ${}^{<\o }2$.
So, by (a), $D$ is a dense set in the poset $\la \L ,\leq _L \ra$.
\hfill $\Box$
\begin{cla}\label{T2653}
For each $n\in \N$, $\p \in {}^{n-1}2$ and $k,j\in 2$ we have
\begin{equation}\label{EQ2722}\textstyle
q^{\bigcup _{i<n_{\p ^{\smallfrown}k} }L_i}_{K_{\p ^{\smallfrown} k {}^{\smallfrown}j}} \sim d_n \Leftrightarrow k=1.
\end{equation}
\end{cla}
\dok
($\Leftarrow$) By (iv) we have $q^{\bigcup _{i<n_{\p ^{\smallfrown}1} }L_i}_{K_{\p ^{\smallfrown} 1 {}^{\smallfrown}j}}
<_L q^{\bigcup _{i<l_n }L_i}_{K_{\p ^{\smallfrown} 1 }\cup \{ d_n\}}$, which, by (\ref{EQ2605}), implies
$d_n \in K_{\p ^{\smallfrown} 1 {}^{\smallfrown}j}$. By (L4),
$q^{\bigcup _{i<n_{\p ^{\smallfrown}1} }L_i}_{K_{\p ^{\smallfrown} 1 {}^{\smallfrown}j}}\in
R^{\bigcup _{i<n_{\p ^{\smallfrown}1} }L_i}_{K_{\p ^{\smallfrown} 1 {}^{\smallfrown}j}}$ so we have
$q^{\bigcup _{i<n_{\p ^{\smallfrown}1} }L_i}_{K_{\p ^{\smallfrown} 1 {}^{\smallfrown}j}} \sim d_n $.

($\Rightarrow$) By (iii) we have
$q^{\bigcup _{i<n_{\p ^{\smallfrown}0} }L_i}_{K_{\p ^{\smallfrown} 0 {}^{\smallfrown}j}}
<_L q^{\bigcup _{i<l_n }L_i}_{K_{\p ^{\smallfrown} 0 }}$, which, by  (\ref{EQ2605}) and (i), implies
\begin{equation}\label{EQ2723}\textstyle
K_{\p ^{\smallfrown} 0 {}^{\smallfrown}j} \cap \bigcup _{i<l_n }L_i =K_{\p ^{\smallfrown} 0 }\subset \bigcup _{i<n_\p }L_i.
\end{equation}
By (\ref{EQ2714}) and (\ref{EQ2719}) we have
\begin{equation}\label{EQ2724}\textstyle
d_n \in L_{n_{\p _n \upharpoonright (n-1)}}\cap \bigcup _{i<l_n }L_i
\end{equation}
and, by (ii), $n_\p \leq n_{\p _n \upharpoonright (n-1)}< l_n <n_{\p ^\smallfrown 0}$, which implies $d_n \not\in \bigcup _{i<n_\p }L_i$
and, by (\ref{EQ2723}) and (\ref{EQ2724}),
$d_n \in \bigcup _{i<l_n }L_i\setminus K_{\p ^{\smallfrown} 0 {}^{\smallfrown}j}
\subset \bigcup _{i< n_{\p ^{\smallfrown} 0}}L_i\setminus K_{\p ^{\smallfrown} 0 {}^{\smallfrown}j}$. Now, since
$q^{\bigcup _{i<n_{\p ^{\smallfrown}0} }L_i}_{K_{\p ^{\smallfrown} 0 {}^{\smallfrown}j}} \in
R^{\bigcup _{i<n_{\p ^{\smallfrown}0} }L_i}_{K_{\p ^{\smallfrown} 0 {}^{\smallfrown}j}}$, we have
$q^{\bigcup _{i<n_{\p ^{\smallfrown}0} }L_i}_{K_{\p ^{\smallfrown} 0 {}^{\smallfrown}j}} \not\sim d_n$.
Claim \ref{T2653} is proved.
\hfill $\Box$
\begin{cla}\label{T2654}
For each $\f \in {}^{<\o }2 \setminus \{ \emptyset \}$, each $j\in 2$, and each $k\in \dom \f$ we have
\begin{equation}\label{EQ2725}\textstyle
q^{\bigcup _{i<n_\f } L_i}_{K_{\f ^{\smallfrown} j }} \sim d_{k+1} \Leftrightarrow \f (k)=1.
\end{equation}
\end{cla}
\dok
Let $k\in \dom \f$. Then there is $s\in 2$ such that $(\f \upharpoonright k) ^{\smallfrown} \f (k) {}^{\smallfrown}s \subset \f ^{\smallfrown} j$
and by (vii) we have
$f(\f ^{\smallfrown} j) \leq _L
f((\f \upharpoonright k )^{\smallfrown} \f (k) {}^{\smallfrown}s)$, that is,
\begin{equation}\label{EQ2726}\textstyle
q^{\bigcup _{i<n_\f } L_i}_{K_{\f ^{\smallfrown} j }}
\leq _L q^{\bigcup _{i<n_{(\f \upharpoonright k )^{\smallfrown} \f (k)} } L_i}_{K_{ (\f \upharpoonright k ) ^{\smallfrown} \f (k) {}^{\smallfrown}s}}
\end{equation}
and, by Claim \ref{T2653},
\begin{equation}\label{EQ2727}\textstyle
 q^{\bigcup _{i<n_{(\f \upharpoonright k )^{\smallfrown} \f (k)} } L_i}_{K_{ (\f \upharpoonright k ) ^{\smallfrown} \f (k) {}^{\smallfrown}s}}
 \sim d_{k+1} \Leftrightarrow \f (k)=1.
\end{equation}
By (ii) we have $n_{\p _{k+1}\upharpoonright k}< l_{k+1}< n_{(\f \upharpoonright k )^{\smallfrown} \f (k)}$
so, by (\ref{EQ2724}), $d_{k+1}\in L_{n_{\p _{k+1}\upharpoonright k}}\subset \bigcup _{i<n_{(\f \upharpoonright k) ^{\smallfrown} \f (k)} } L_i$
and, by (\ref{EQ2726}) and Fact \ref{T2627},
$q^{\bigcup _{i<n_\f } L_i}_{K_{\f ^{\smallfrown} j }}\sim d_{k+1}$ if and only if
$q^{\bigcup _{i<n_{(\f \upharpoonright k )^{\smallfrown} \f (k)} } L_i}_{K_{ (\f \upharpoonright k ) ^{\smallfrown} \f (k) {}^{\smallfrown}s}}\sim d_{k+1}$.
Now (\ref{EQ2725}) follows from (\ref{EQ2727}) and Claim \ref{T2654} is proved.
\hfill $\Box$
\begin{cla}\label{T2662}
$D\in \P (R)$.
\end{cla}
\dok
We show that the subgraph $D =\{ d_n: n\in \N \}$ of $R$ satisfies condition (\ref{EQ2740}), that is, for $K\subset H \in [\N ]^{<\o}$ we show that
$D\cap R^{\{ d_n : n\in H \}}_{\{ d_n : n\in K \}}\neq \emptyset$. Let $m:= \max H$ and let $\f \in {}^m 2$ be defined by
\begin{equation}\label{EQ2742}\textstyle
\f (k)=
\left\{ \begin{array}{ll}
              1                             & \mbox{if } k+1 \in K,\\
              0                               & \mbox{if } k+1 \in [1,m]\setminus K.

        \end{array}
\right.
\end{equation}
Then, by Claim \ref{T2654}, for $\;l<m$ we have
\begin{equation}\label{EQ2741}\textstyle
q:= q^{\bigcup _{i< n_\f }L_i}_{K_{\f ^{\smallfrown}0 }} \sim d_{l+1}\Leftrightarrow \f (l)=1 .
\end{equation}
Thus for $l\in K$ by (\ref{EQ2742}) we have $\f (l-1)=1$  and, by (\ref{EQ2741}), $q\sim d_l$;
similarly, for $l\in H\setminus K$ by (\ref{EQ2742}) we have $\f (l-1)=0$  and, by (\ref{EQ2741}), $q\not\sim d_l$. So,
$q\in R^{\{ d_n : n\in H \}}_{\{ d_n : n\in K \}}$.
Since $q\in \L$ and, by Claim \ref{T2661}(b), $D$ is a dense set in the poset $\la \L ,\leq _L \ra$, there is
$d_r\in D$ such that $d_r \leq _L q^{\bigcup _{i< n_\f }L_i}_{K_{\f ^{\smallfrown}0 }}$.
By Fact \ref{T2627} for each $u\in\bigcup _{i< n_\f }L_i$ we have:
$d_r\sim u$ if and only if
$q\sim u$; so it remains to be shown that
$\{ d_n : n\in H \} \subset \bigcup _{i< n_\f }L_i$.
But this is true since by (ii) we have $l_m<n_\f$ and for $n\in H$ we have $n\leq m$ so, by (\ref{EQ2714}),
$d_n \in \bigcup _{i< l_n }L_i\subset \bigcup _{i< l_m }L_i\subset \bigcup _{i< n_\f }L_i$.
Thus $d_r\in D\cap R^{\{ d_n : n\in H \}}_{\{ d_n : n\in K \}}$ and Claim \ref{T2662} is proved.
\hfill $\Box$
\begin{cla}\label{T2655}
For each $\f \in {}^{<\o }2$ and each $j\in 2$  we have
\begin{equation}\label{EQ2728}\textstyle
\Big(-\infty , q^{\bigcup _{i< n_\f }L_i}_{K_{\f ^{\smallfrown}j }} \Big)_{\la L , \leq _L \ra} \cap D \in \P (R).
\end{equation}
\end{cla}
\dok
Let $K\subset H $ be finite subsets of $(-\infty , q^{\bigcup _{i< n_\f }L_i}_{K_{\f ^{\smallfrown}j }} )_{\la L , \leq _L \ra} \cap D$, where
$K=\{ d_k : k\in F\}$ and $H\setminus K=\{ d_k : k\in G\}$. Then $F$ and $G$ are disjoint finite subsets of $\N$. By (vii), for $\th \in {}^{<\o }2$ and $r\in 2$
we have
\begin{equation}\label{EQ2729}\textstyle
\f ^{\smallfrown} j \varsubsetneq \th ^{\smallfrown} r
\Leftrightarrow
 q^{\bigcup _{i< n_\th }L_i}_{K_{\th ^{\smallfrown}r }}
< _L q^{\bigcup _{i< n_\f }L_i}_{K_{\f ^{\smallfrown}j }}.
\end{equation}
By the assumption and (\ref{EQ2719}), for $k\in F\cup G$ we have
$d_k = q^{\bigcup _{i<n_{\p _k\upharpoonright (k-1)}}L_i}_{K_{(\p _k\upharpoonright (k-1))^{\smallfrown}0}}
<_L q^{\bigcup _{i< n_\f }L_i}_{K_{\f ^{\smallfrown}j }}$, which by (\ref{EQ2729}) implies
$\f ^{\smallfrown} j \varsubsetneq (\p _k\upharpoonright (k-1))^{\smallfrown}0$ and, hence, $|\f |+1 <k$.
Thus $k-1\not \in \dom (\f ^{\smallfrown}j)$, for each $k\in F\cup G$ and we take $\th \in {}^{< \o}2$ such that $H\subset \bigcup _{i<n_\th}L_i$ and
\begin{equation}\label{EQ2730}\textstyle
\f ^{\smallfrown}j \subset \th \;\; \land  \;\; \forall k\in F \; (\th (k-1)=1 )  \;\; \land \;\;  \forall k\in G \; (\th (k-1)=0 ).
\end{equation}
Let $r\in 2$ .
By Claim \ref{T2654} and (\ref{EQ2730}), for each $d_k\in H$ we have $q^{\bigcup _{i< n_\th }L_i}_{K_{\th ^{\smallfrown}r }} \sim d_k$ iff $\th (k-1)=1$ iff $k\in F$ iff $d_k\in K$; thus
\begin{equation}\label{EQ2731}\textstyle
q^{\bigcup _{i< n_\th }L_i}_{K_{\th ^{\smallfrown}r }}\in L^{\bigcup _{i< n_\th }L_i}_{K_{\th ^{\smallfrown}r }}\subset R^H_K .
\end{equation}
By (\ref{EQ2730}) and (\ref{EQ2729}) we have  $q^{\bigcup _{i< n_\th }L_i}_{K_{\th ^{\smallfrown}r }}
< _L q^{\bigcup _{i< n_\f }L_i}_{K_{\f ^{\smallfrown}j }}$ and, since $D$ is dense in $\la \L \leq _L\ra$, there is
$d_{k_0}\leq_L q^{\bigcup _{i< n_\th }L_i}_{K_{\th ^{\smallfrown}r }}$
which, by (\ref{EQ2731}) and Fact \ref{T2627}, implies $d_{k_0}\in R^H_K$.
Clearly, $d_{k_0}\in (-\infty , q^{\bigcup _{i< n_\f }L_i}_{K_{\f ^{\smallfrown}j }} )_{\la L , \leq _L \ra} \cap D$ and Claim \ref{T2655} is proved.
\hfill $\Box$
\begin{cla}\label{T2656}
For each $n\in \N$ the set $\A _n =\{ (-\infty , q)_{\la L, \leq _L \ra} \cap D :q\in \L _n \}$, that is
$$
\A _n = \Big\{ (-\infty , q^{\bigcup _{i< n_\f }L_i}_{K_{\f ^{\smallfrown}j }} )_{\la L , \leq _L \ra} \cap D : \f \in {}^{n-1 }2 \land j\in 2\Big\}
$$
is an antichain in $\la \P (R) , \subset \ra$, maximal below $D$.
\end{cla}
\dok
By (vi) and (vii), $\L _n$ is a maximal antichain in the reversed tree $\la \L , \leq _L \ra$ of size $2^n$ and, hence, the sets
$(-\infty , q)_{\la L, \leq _L \ra}$, $q\in \L _n$, are disjoint and, by Claim \ref{T2655}, $\A _n$ is an antichain in $\la \P (R) , \subset \ra$.
Since $\L \setminus \bigcup _{q\in \L _n}(-\infty , q)_{\la L, \leq _L \ra}\subset \bigcup _{i\leq n}\L _n$ is a finite set,
for $E\in \P (D)$ we have $E_1=E\setminus \bigcup _{i\leq n}\L _n \in \P (D)$ (see Fact \ref{T2600}(a)) and, since $E_1\subset\bigcup \A _n$,
by Fact \ref{T2600}(b)
$E_1$ is compatible with some element of $\A _n$. Claim \ref{T2656} is proved.
\hfill $\Box$
\begin{cla}\label{T2657}
The mapping $F : \la\L ,\leq _L \ra \rightarrow \la {}^{<\o }\o , \supset\ra$ defined by $F(q^\emptyset _\emptyset )=\emptyset$ and
\begin{equation}\label{EQ2732}\textstyle
F(q^{\bigcup _{i< n_\f }L_i}_{K_{\f ^{\smallfrown}j }})=\Big\la m ^k_{K_{\f ^{\smallfrown}j}\cap \bigcup _{i<k}L_i}: k\in \{0,\dots ,n_\f\} \Big\ra ,
\end{equation}
for $\f \in {}^{<\o }2$ and $j\in 2$, is an embedding.
\end{cla}
\dok
By (i) we have $K_{\f ^{\smallfrown}j }\subset \bigcup _{i< n_\f }L_i$ thus $K_{\f ^{\smallfrown}j }\cap \bigcup _{i<k}L_i \subset \bigcup _{i< k }L_i$, for
all $k\leq n_\f$,
and, by (\ref{EQ2703}), $m ^k_{K_{\f ^{\smallfrown}j}\cap \bigcup _{i<k}L_i}\in \o$.
Thus $F$ is well defined.

For a proof that $F$ is an injection we suppose that
 $q^{\bigcup _{i< n_\f }L_i}_{K_{\f ^{\smallfrown}j }}\neq q^{\bigcup _{i< n_\p }L_i}_{K_{\p ^{\smallfrown}l }}$.

If $n_\f \neq n_\p$ then the corresponding sequences are of different length and, hence, different.

Let $n_\f = n_\p$. Suppose that $|\f |\neq |\p |$, say $|\f |< |\p |$. Then, by (ii), $n_\f <l_{|\f |+1}\leq l_{|\p |}<n_\p$, which contradicts our assumption.
Thus $\f ,\p \in {}^{n-1}2$, for some $n\in \N$.

If $\f = \p$, then, for example $j=0$ and $l=1$ and, by (v),
$m^{n_\f }_{K_{\f ^{\smallfrown} 0}}\neq m^{n_\f }_{K_{\f ^{\smallfrown} 1}}$ so the corresponding sequences
have the last elements different.

If $\f \neq \p$, then there are $k<n-1$ and $\e\in {}^{k-1}2$ such that, for example,
$\e ^{\smallfrown} 0 \subset \f$ and $\e ^{\smallfrown} 1 \subset \p$.
Since $\e ^{\smallfrown} 0 \subset \f ^{\smallfrown} j$ and $\e ^{\smallfrown} 1 \subset \p ^{\smallfrown} l$ by (vii) we have
$q^{\bigcup _{i< n_\f }L_i}_{K_{\f ^{\smallfrown}j }}\leq _L q^{\bigcup _{i< n_\e }L_i}_{K_{\e ^{\smallfrown}0 }} $
and $q^{\bigcup _{i< n_\p }L_i}_{K_{\p ^{\smallfrown}l }}\leq _L q^{\bigcup _{i< n_\e }L_i}_{K_{\e ^{\smallfrown}1 }} $,
which by (\ref{EQ2605}) implies
\begin{equation}\label{EQ2734}\textstyle
K_{\f ^{\smallfrown}j }\cap \bigcup _{i< n_\e }L_i=K_{\e ^{\smallfrown}0} \;\mbox{ and } \;\;
K_{\p ^{\smallfrown}l }\cap \bigcup _{i< n_\e }L_i=K_{\e ^{\smallfrown}1}.
\end{equation}
By (v) we have
$m^{n_\e }_{K_{\e ^{\smallfrown} 0 }}\neq m^{n_\e }_{K_{\e ^{\smallfrown} 1 }}$ so, by (\ref{EQ2734}),
$m^{n_\e }_{K_{\f ^{\smallfrown}j }\cap \bigcup _{i< n_\e }L_i }\neq m^{n_\e }_{ K_{\p ^{\smallfrown}l }\cap \bigcup _{i< n_\e }L_i}$,
which shows that $F(q^{\bigcup _{i< n_\f }L_i}_{K_{\f ^{\smallfrown}j }})\neq F(q^{\bigcup _{i< n_\p }L_i}_{K_{\p ^{\smallfrown}l }})$.
So, $F$ is an injection indeed.

For a proof that $F$ is order preserving we suppose that $q^{\bigcup _{i< n_\f }L_i}_{K_{\f ^{\smallfrown}j }}\leq _L q^{\bigcup _{i< n_\p }L_i}_{K_{\p ^{\smallfrown}l }}$. Then $n_\f \geq n_\p$
and $K_{\f ^{\smallfrown}j }\cap \bigcup _{i< n_\p }L_i=K_{\p ^{\smallfrown}l }$ so, for each $k\leq n_\p$ we have
$K_{\f ^{\smallfrown}j }\cap \bigcup _{i< k }L_i
=K_{\f ^{\smallfrown}j }\cap \bigcup _{i< n_\p }L_i\cap \bigcup _{i< k }L_i
=K_{\p ^{\smallfrown}l }\cap \bigcup _{i< k }L_i$ which implies
\begin{equation}\label{EQ2735}\textstyle
\Big\la m ^k_{K_{\p ^{\smallfrown}l}\cap \bigcup _{i<k}L_i}: k\leq n_\p \Big\ra \subset
\Big\la m ^k_{K_{\f ^{\smallfrown}j}\cap \bigcup _{i<k}L_i}: k\leq n_\f \Big\ra
\end{equation}
that is $F(q^{\bigcup _{i< n_\f }L_i}_{K_{\f ^{\smallfrown}j }})\supset F( q^{\bigcup _{i< n_\p }L_i}_{K_{\p ^{\smallfrown}l }})$.
Thus the mapping $F$ is order preserving.

Suppose that (\ref{EQ2735}) holds and that $q^{\bigcup _{i< n_\f }L_i}_{K_{\f ^{\smallfrown}j }}\not\leq _L q^{\bigcup _{i< n_\p }L_i}_{K_{\p ^{\smallfrown}l }}$.
Then, since $F$ is one-to-one, $q^{\bigcup _{i< n_\p }L_i}_{K_{\p ^{\smallfrown}l }}\not\leq _L q^{\bigcup _{i< n_\f }L_i}_{K_{\f ^{\smallfrown}j }}$
and, by (vii), $\f ^{\smallfrown}j$ and $\p ^{\smallfrown}l $ are incompatible elements of ${}^{<\o}2$.
Hence there is $\e\in {}^{<\o}2$ such that, for example,
$\e ^{\smallfrown} 0 \subset \f^{\smallfrown} j$ and $\e ^{\smallfrown} 1 \subset \p ^{\smallfrown} l$
and, like in the proof that $F$ is an injection we obtain
$m^{n_\e }_{K_{\f ^{\smallfrown}j }\cap \bigcup _{i< n_\e }L_i }\neq m^{n_\e }_{ K_{\p ^{\smallfrown}l }\cap \bigcup _{i< n_\e }L_i}$,
which contradicts (\ref{EQ2735}) (by (ii), $\e \subset \p$ implies $n_\e \leq n_\p$). Claim \ref{T2657} is proved.
\kdok

\noindent
By Claims \ref{T2661}(a) and \ref{T2657} the composition $g= F\circ f :\la {}^{<\o }2 ,\supset \ra \rightarrow  \la {}^{<\o }\o , \supset\ra$, where
$$
\la {}^{<\o }2 ,\supset \ra \cong _f  \la\L ,\leq _L \ra \hookrightarrow _F \la {}^{<\o }\o , \supset\ra
$$
is an embedding and by  Lemma \ref{T2639} we have $\T :=g[{}^{<\o}2]\!\uparrow \;\in \Bt ({}^{<\o}\o )$.
Clearly,
$$
\T =\Big\{ \p \in {}^{<\o}\o :
\exists \f \in {}^{<\o }2 \; \exists j\in 2 \;\; \p \subset \la m ^k_{K_{\f ^{\smallfrown}j}\cap \bigcup _{i<k}L_i}: k\leq n_\f \ra\Big\}.
$$
\begin{cla}\label{T2658}
$D\Vdash \t \upharpoonright \check{n} \in \check{\T }$, for each $n\in \o$.
\end{cla}
\dok
Let $n\in \o$. We have to prove that
\begin{equation}\label{EQ2736}\textstyle
D\Vdash \exists \f \in \check{{}^{<\o }2} \;\; \exists j\in \check{2} \;\;
\Big(\check{n}\leq n_\f \;\land\; \forall k<\check{n} \;\;\t (k)= \check{m ^k_{K_{\f ^{\smallfrown}j}\cap \bigcup _{i<k}L_i}}\Big).
\end{equation}
Let $G$ be a generic filter and $D\in G$. By (ii), there is $n_0\in \N$ such that $n<l_{n_0 -1}$ and, by (ii) again,
\begin{equation}\label{EQ2737}\textstyle
\forall \f \in {}^{n_0 -1}2 \;\; n_\f > l_{n_0 -1} >n .
\end{equation}
By Claim \ref{T2656}, the set
$\A _{n_0} = \{ (-\infty , q^{\bigcup _{i< n_\f }L_i}_{K_{\f ^{\smallfrown}j }} )_{\la L , \leq _L \ra} \cap D : \f \in {}^{n_0 -1 }2 \land j\in 2\}$
is an antichain in $\la \P (R),\subset \ra$, maximal below $D$ and, since $D\in G$, there are $\f _0 \in {}^{n_0 -1 }2$ and $j_0\in 2$ such that
\begin{equation}\label{EQ2738}\textstyle
\Big( -\infty , q^{\bigcup _{i< n_{\f _0} }L_i}_{K_{\f _0 ^{\smallfrown}j_0 }} \Big)_{\la L , \leq _L \ra} \cap D \in G.
\end{equation}
By (\ref{EQ2737}), for $k<n$ we have $k<n_{\f _0}$ and, hence,
$q^{\bigcup _{i< n_{\f _0} }L_i}_{K_{\f _0 ^{\smallfrown}j_0 }}\leq _L q^{\bigcup _{i< k }L_i}_{K_{\f _0 ^{\smallfrown}j_0 }\cap \bigcup _{i< k }L_i }$.
So, by Fact \ref{T2627},
$(-\infty , q^{\bigcup _{i< n_{\f _0} }L_i}_{K_{\f _0 ^{\smallfrown}j_0 }} )_{\la L , \leq _L \ra} \cap D
\subset (-\infty , q^{\bigcup _{i< k }L_i}_{K_{\f _0 ^{\smallfrown}j_0 }\cap \bigcup _{i< k }L_i } ]_{\la L , \leq _L \ra}$
$=L^{\bigcup _{i< k }L_i}_{K_{\f _0 ^{\smallfrown}j_0 }\cap \bigcup _{i< k }L_i }$
and, by (\ref{EQ2703}),
$L^{\bigcup _{i< k }L_i}_{K_{\f _0 ^{\smallfrown}j_0 }\cap \bigcup _{i< k }L_i }
\Vdash \t (\check{k})= (m ^k_{K_{\f _0 ^{\smallfrown}j_0}\cap \bigcup _{i<k}L_i})\check{\;}$,
which, together with (\ref{EQ2738}), gives
$\t _G (k)= m ^k_{K_{\f _0 ^{\smallfrown}j_0}\cap \bigcup _{i<k}L_i}$.
(\ref{EQ2736}) is proved.
\kdok
Now, since $D\subset \L \subset L \subset B$, (\ref{EQ2702}) follows from Claims \ref{T2662} and  \ref{T2658}.
\kdok
\section{Appendix}
In this section we establish the mentioned forcing-free translation of the 2-localization property to the language of Boolean algebras.
We recall that a complete Boolean algebra $\B$ is said to have
the 2-localization property iff
\begin{equation}\label{EQ2746}\textstyle
1_\B \Vdash _\B \forall x:\check{\o} \rightarrow \check{\o } \;\;
                \exists \T \in  ((\Bt ({}^{<\o }\o))^V )\check{\;}\;\;
                \forall n\in \check{\o } \;\; x\upharpoonright n\in \T .
\end{equation}
\begin{te}\label{T2664}
A complete Boolean algebra $\B$ has the 2-localization property if and only if for each double sequence $[b_{n m }: \langle n , m \rangle \in \o \times \o ]$ of
elements of ${\mathbb B}$ we have
\begin{equation}\label{EQ2747}
\textstyle
\bigwedge _{n \in \o }\;
\bigvee _{m \in \o }\;
b_{n m}
=
\bigvee _{\T \,\in\, \Bt ({}^{<\o }\o)}\;
\bigwedge _{n \in \o }\;
\bigvee _{\f \,\in\, \T \cap {}^{n+1}\o }\;
\bigwedge _{k\leq n}\;
b_{k\f (k)}.
\end{equation}
\end{te}
\dok
For notational simplicity, let $b$ and $c$ denote the left hand side and the right hand side of (\ref{EQ2747}) respectively.

First we prove that $b\geq c$
holds for each sequence $[b_{nm}]$ in any complete Boolean algebra
(and, hence, that the equality (\ref{EQ2747}) is equivalent to the inequality $b\leq c$).
For $\T \in \Bt ({}^{<\o }\o)$ and $n\in \o$ we have: $\bigwedge _{k\leq n}b_{k\f (k)}\leq b_{n\f (n)}$, for each $\f \,\in\, \T \cap {}^{n+1}\o $. Thus
$\bigvee _{\f \,\in \, \T \cap {}^{n+1}\o }\;
\bigwedge _{k\leq n}\;
b_{k\f (k)}
\leq
\bigvee _{\f \,\in \, \T \cap {}^{n+1}\o }\;
b_{n\f (n)}
\leq
\bigvee _{m \in \o }\;
b_{n m}
$.
This holds for all $n\in \o $ and, hence,
$$\textstyle
\bigwedge _{n \in \o }\;
\bigvee _{\f \,\in\, \T \cap {}^{n+1}\o }\;
\bigwedge _{k\leq n}\;
b_{k\f (k)}
\leq
\bigwedge _{n \in \o }\;
\bigvee _{m \in \o }\;
b_{n m}
= b .
$$
This holds for all $\T \in \Bt ({}^{<\o }\o)$ so we have $c\leq b$. Now we prove the theorem.

($\Leftarrow$) Assuming that (\ref{EQ2747}) holds for each double sequence $[b_{nm}]$ in $\B$ we prove that $\B$ has
the $2$-localization property. Let $G$ be a $\B$-generic filter over $V$ and $x\in V_\B [G]$, where $x:\o \rightarrow \o$.
Let $\t$ be a $\B$-name and $b_0\in G$, where $x=\t _G$ and $b_0 \Vdash \t : \check{\o }\rightarrow \check{\o }$. Defining
$b_{nm}:=\| \t (\check{n})= \check{m}\|=\| \la \check{n},\check{m} \ra {} \in \t\|$ we obtain
$b_0 \leq \| \forall n\in \check{\o }\; \exists m\in \check{\o }\;  \la n,m \ra\in \t \|= \bigwedge _{n \in \o }\;
\bigvee _{m \in \o }\;
b_{n m} =c$, which means that
$b_0 \Vdash
\exists \T \,\in\, ((\Bt ({}^{<\o }\o))^ V)\check{\;}\;
\forall n \in \check{\o } \;
\exists \f \,\in\, \T \cap {}^{n+1}\o \;
\forall k\leq n\;
\t (k )=\f (k)$.
So, since $b_0\in G$, there is $\T \,\in\, \Bt ({}^{<\o }\o)\cap V$ such that for each $n\in \o$ we have $\t \upharpoonright (n+1)\in \T$.
Clearly, $\t \upharpoonright 0\in \T$ and we are done.

($\Rightarrow$) Suppose that the algebra $\B$ has the 2-localization property. For a double sequence $[b_{nm}]$ in $\B$ we prove (\ref{EQ2747}), that is,
$b\leq c$. If $b=0$, we are done. Otherwise, defining the sequence  $[b_{nm}^*]$ by $b_{n0}^*=b_{n0}$ and $b_{nm}^*=b_{nm}\setminus \bigvee _{i<m}b_{ni}$, for
$m>0$, for each $n\in \o$ we have
$$\textstyle
b\leq \bigvee _{m\in \o}b_{nm}^*
$$
$$\textstyle
m_1 \neq m_2 \Rightarrow b_{n m_1}^* \wedge b_{n m_2}^*=0.
$$
Hence for a $\B$-name
$$
\t := \Big\{ \Big\la \la n,m \ra \check {\;} ,b_{n m}^* \Big\ra : \la n,m \ra \in \o \times \o \Big\}
$$
we have $b\Vdash \t : \check{\o}\rightarrow \check{\o} $ and $\| \t (\check{n})= \check{m}\|=b_{n m}^*$, which, by the assumption,
implies that
$\textstyle
b \Vdash _\P \exists \T \in  ((\Bt ({}^{<\o }\o))^V )\check{\;}\;\;
                \forall n\in \check{\o } \;\; x\upharpoonright (n+1)\in \T .
$
Thus
$$\textstyle
b\leq \bigvee _{\T \,\in\, \Bt ({}^{<\o }\o)}\;
\bigwedge _{n \in \o }\;
\bigvee _{\f \,\in\, \T \cap {}^{n+1}\o }\;
\bigwedge _{k\leq n}\;
b_{k\f (k)}^* \leq c,
$$
because $b_{nm}^* \leq b_{nm}$, for all $n,m\in\o$, and the proof is over.
\kdok

\vspace{2mm}
\noindent
{\bf Acknowledgement}

\vspace{2mm}
Research of M.\  Kurili\'c was supported by the Ministry of Education and Science of the Republic of Serbia (Project 174006).

Research of S.\ Todor\v cevi\'c  was supported by the grants from CNRS and NSERC.


\footnotesize
\begin{thebibliography}{99}
\bibitem{Camer}
      P.\ J.\ Cameron,
      The random graph,
      The mathematics of Paul Erd\" os II,
      Algorithms Combin. 14 (Springer, Berlin, 1997) 333--351.
\bibitem{Erdos2}
      P.\ Erd\H{o}s, A.\ R\'enyi,
      Asymmetric graphs,
      Acta Math.\ Acad.\ Sci.\ Hungar., 14 (1963) 295--315.
\bibitem{Jech}
      T.\ Jech,
      Set theory. Second edition.
      Perspectives in Mathematical Logic. Springer-Verlag, Berlin, 1997.
\bibitem{Ktow}
      M.\ S.\ Kurili\'c,
      From $A_1$ to $D_5$: Towards a forcing-related classification of relational structures,
      J.\ Symbolic Logic, 79,1 (2014) 279--295.
\bibitem{Kur1}
      M.\ S.\ Kurili\'c,
      Posets of copies of countable scattered linear orders,
      Ann.\ Pure  Appl.\ Logic, 165 (2014)  895--912.
\bibitem{Kurord}
      M.\ S.\ Kurili\'c,
      Forcing with copies of countable ordinals,
      Proc.\ Amer.\ Math.\ Soc.\ to appear.\\
      http://arxiv.org/abs/1304.7714
\bibitem{Kurstr}
      M.\ S.\ Kurili\'c,
      Isomorphic and strongly connected components,
      Arch.\ Math.\ Logic, to appear. DOI10.1007/s00153-014-0399-2\\
      http://arxiv.org/abs/1311.5049
\bibitem{Kuremb}
      M.\ S.\ Kurili\'c,
      Maximally embeddable components,
      Arch.\ Math.\ Logic 52,7 (2013) 793-808.
\bibitem{Kdif}
      M.\ S.\ Kurili\'c,
      Different Similarities, to appear.
\bibitem{Kmon}
      M.\ S.\ Kurili\'c,
      Monoids of self-embeddings, to appear.
\bibitem{KurTod}
      M.\ S.\ Kurili\'c, S.\ Todor\v cevi\'c,
      Forcing by non-scattered sets,
      Ann.\ Pure Appl.\ Logic 163 (2012) 1299--1308.
\bibitem{KurTodR1}
      M.\ S.\ Kurili\'c, S.\ Todor\v cevi\'c,
      Copies of the random graph, to appear.\\
      http://arxiv.org/abs/1410.6320
\bibitem{NewRos}
      L.\ Newelski, A.\ Ros{\l}anowski,
      The ideal determined by the unsymmetric game,
      Proc.\ Amer.\ Math.\ Soc., 117,3 (1993) 823--831.
\bibitem{Rado}
      R.\ Rado,
      Universal graphs and universal functions,
      Acta Arith., 9 (1964) 331--340.
\end{thebibliography}
\end{document}